\documentclass[twoside,12pt]{article}
\usepackage{amssymb}
\usepackage{amsthm}
\usepackage{amsmath}
\usepackage{a4wide}
\usepackage[all]{xy}
\usepackage{enumerate}
\usepackage{fancyhdr}

\newtheorem{theorem}{Theorem}[section]
\newtheorem{proposition}[theorem]{Proposition}
\newtheorem{corollary}[theorem]{Corollary}
\newtheorem{lemma}[theorem]{Lemma}
\theoremstyle{definition}

\newtheorem*{Beweis}{Proof}
\newtheorem{definition}[theorem]{Definition}
\newtheorem{punto}[theorem]{}
\theoremstyle{remark}
\newtheorem{remark}[theorem]{Remark}

\newtheorem{remarks}[theorem]{Remarks}
\swapnumbers
\input{tcilatex}
\begin{document}

\title{A Note on Coseparable Coalgebras\textbf{\thanks{%
MSC (2000): 16W30 \newline Keywords: Coseparable Coalgebras; Coseparable Corings; Separable Algebras; Separable Rings; Separable Functors}}}
\author{\textbf{Jawad Y. Abuhlail\thanks{%
Supported by King Fahd University of Petroleum $\&\ $Minerals (KFUPM)}} \\
Department of Mathematics and Statistic, Box $\#\ $5046\\
King Fahd University of Petroleum $\&$ Minerals\\
31261 Dhahran - Saudi Arabia\\
abuhlail@kfupm.edu.sa}
\date{}
\maketitle

\begin{abstract}
Given a coalgebra $C$ over a commutative ring $R,$ we show that $C$ can be
considered as a (\emph{not necessarily counital}) $C^{\ast op}$-coring.
Moreover, we show that this coring has a left (right) counity if and only if
$C$ is coseparable as an $R$-coalgebra.
\end{abstract}

\section{Introduction}

\qquad \emph{Coseparable coalgebras} (\emph{coseparable corings}) are dual
to separable algebras (separable rings), and were introduced by Larson \cite%
{Lar1973} (Guzm\'{a}n \cite{Gut1989}). These were investigated, using \emph{%
homological, categorical and module theoretical }approaches by several
authors including Doi \cite{Doi1981}, Casta\~{n}o Iglesias et. al. \cite%
{C-IG-TN1997} and Brzezi\'{n}ski et. al. \cite{BTW2006} (G\'{o}%
mez-Torrecillas et. al. in \cite{G-TL2003}, \cite{G-T-2002} and Brzezi\'{n}%
ski et. al. \cite{BKW-2005}). Generalizations were also considered by
Nakajima (e.g. \cite{Nak1979}, \cite{Nak1980}), who introduced the notion of%
\emph{\ coseparable coalgebras over coalgebras}.

The main goal of this short note is to recover a coseparable coalgebra $C$
as left (right) counital coring over the opposite dual algebra $C^{\bullet
}:=C^{\ast op}.$ After this brief introduction, we provide in Section two
some preliminaries about (not necessarily counital) corings and their
categories of comodules. We introduce also the notions of left (right)
counital rings and their dual right (left) unital rings. In Section three,
we show that every coalgebra $C$ can be considered (in a natural way) as a
(not necessarily counital) coring over $C^{\bullet }:=C^{\ast op}.$ In
Proposition \ref{C-separable}, we show that an $R$-coalgebra $C,$ for which
the associated coring $(C:C^{\bullet })$ is left (right) counital is a
separable $R$-algebra as well as a coseparable $R$-algebra and induce a
\emph{separable measuring pairing}. In Section four, we prove the main
result, Theorem \ref{cosep-counit}, which states that the coring $%
(C:C^{\bullet })$ is left (right) counital if and only if $C$ is a
coseparable $R$-coalgebra.

\section{Preliminaries}

\qquad Throughout, $R$ denotes a \emph{commutative }ring with $1_{R}\neq
0_{R}$ and $A,B$ are associative (not necessarily commutative) unital $R$%
-algebras. All modules over unital rings are assumed to be unitary. Let $M$
be a right $B$-module and $N$ a left $B$-module. Then we have canonical
isomorphisms%
\begin{equation*}
M\otimes _{B}B\overset{\vartheta _{M}^{r}}{\simeq }M\text{ and }B\otimes
_{B}N\overset{\vartheta _{N}^{l}}{\simeq }N.
\end{equation*}%
If $\gamma :A\rightarrow B$ is a morphism of $R$-algebras, then we consider $%
M$ as a right $A$-module, $N$ as a left $A$-module with actions induced
through $\gamma $ by those of $B,$ and we have an $R$-linear morphism%
\begin{equation}
\chi _{(M,N)}^{(A,B)}:M\otimes _{A}N\rightarrow M\otimes _{B}N.
\label{chi-MN}
\end{equation}

For a left $A$-module (resp. right $A$-module, $(A,A)$-bimodule) $U,$ we
consider the set $^{\ast }U:=\mathrm{Hom}_{A-}(U,A)$ (resp. $U^{\ast }:=%
\mathrm{Hom}_{-A}(U,A),$ $^{\ast }U^{\ast }:=\mathrm{Hom}_{(A,A)}(U,A)$) of
left $A$-linear (resp. right $A$-linear, $(A,A)$-bilinear) maps from $U$ to $%
A$ as a right $A$-module (resp. a left $A$-module, $(A,A)$-bimodule) in the
canonical way. With \emph{locally projective }modules, we mean those in the
sense of \cite{Z-H-1976}.

\begin{punto}
A (\emph{not necessarily unital})\emph{\ }$A$\textbf{-ring} $(T,\mu _{T}:A)$
is an $(A,A)$-bimodule, with an $(A,A)$-bilinear map $\mu _{T}:T\otimes
_{A}T\rightarrow T,$ called \textbf{multiplication}, such that%
\begin{equation*}
\mu _{T}\circ (\mu _{T}\otimes _{A}\mathrm{id}_{T})=\mu _{T}\circ (\mathrm{id%
}_{T}\otimes _{A}\mu _{T}).
\end{equation*}%
An $A$-ring $T$ is said to be \textbf{left unital} (resp. \textbf{right
unital}, \textbf{unital}), iff there exists a right $A$-linear (resp. left $%
A $-linear, $(A,A)$-bilinear) map, called \textbf{left unity map }(resp.
\textbf{right unity map, unity map})\textbf{\ }$\eta _{T}^{l}:A\rightarrow T$
(resp. $\eta _{T}^{r}:A\rightarrow T,$ $\eta _{T}:A\rightarrow T$), such
that $\mu _{T}\circ (\eta _{T}^{l}\otimes _{A}\mathrm{id}_{T})\equiv
\vartheta _{T}^{l}$ (resp. $\mu _{T}\circ (\mathrm{id}_{T}\otimes _{A}\eta
_{T}^{r})\equiv \vartheta _{T}^{r}$ $,$ $\mu _{T}\circ (\eta _{T}\otimes _{A}%
\mathrm{id}_{T})\equiv \vartheta _{T}^{l}$ and $\mu _{T}\circ (\mathrm{id}%
_{T}\otimes _{A}\eta _{T})\equiv \vartheta _{T}^{r}$).
\end{punto}

For an $A$-ring $T$ and a left (right) $T$-module $M,$ we denote with $%
\sigma \lbrack M]$ Wisbauer's category of $T$\emph{-subgenerated} left
(right) $T$-modules, i.e. the category of submodules of $M$-generated left
(right) $T$-modules (e.g. \cite{Wis-1991} and \cite{Wis1996}).

\section*{Left Counital (Right Counital) Corings}

\begin{punto}
With a (not necessarily counital) $A$\textbf{-coring }$(\mathcal{C},\Delta _{%
\mathcal{C}}:A),$ we mean an $(A,A)$-bimodule $\mathcal{C}$ with an $(A,A)$%
-bilinear map $\Delta _{\mathcal{C}}:\mathcal{C}\rightarrow \mathcal{C}%
\otimes _{A}\mathcal{C},$ called \textbf{comultiplication}, such that%
\begin{equation*}
(\Delta _{\mathcal{C}}\otimes _{A}\mathrm{id}_{\mathcal{C}})\circ \Delta _{%
\mathcal{C}}=(\mathrm{id}_{\mathcal{C}}\otimes _{A}\Delta _{\mathcal{C}%
})\circ \Delta _{\mathcal{C}}.
\end{equation*}
\end{punto}

\begin{punto}
Let $(\mathcal{C},\Delta _{\mathcal{C}})$ be an $A$-coring. We say that $(%
\mathcal{C},\Delta _{\mathcal{C}})\ $is

\textbf{left counital, }iff there exists a \emph{right }$A$\emph{-linear map}
$\varepsilon _{\mathcal{C}}^{l}:\mathcal{C}\rightarrow A$ (called \textbf{%
left counity}) such that%
\begin{equation*}
(\varepsilon _{\mathcal{C}}^{l}\otimes _{A}\mathrm{id}_{\mathcal{C}})\circ
\Delta _{\mathcal{C}}\equiv \vartheta _{\mathcal{C}}^{l},\text{ i.e. }\dsum
\varepsilon _{\mathcal{C}}^{l}(c_{1})c_{2}=c\text{ for all }c\in \mathcal{C};
\end{equation*}

\textbf{right counital, }iff there exists a \emph{left }$A$\emph{-linear map}
$\varepsilon _{\mathcal{C}}^{r}:\mathcal{C}\rightarrow A$ (called \textbf{%
right counity}), such that%
\begin{equation*}
(\mathrm{id}_{\mathcal{C}}\otimes _{A}\varepsilon _{\mathcal{C}}^{r})\circ
\Delta _{\mathcal{C}}\equiv \vartheta _{\mathcal{C}}^{r},\text{ i.e. }\dsum
c_{1}\varepsilon _{\mathcal{C}}^{r}(c_{2})=c\text{ for all }c\in \mathcal{C};
\end{equation*}

\textbf{counital}, iff there exists an $(A,A)$\emph{-bilinear map} $%
\varepsilon _{\mathcal{C}}:\mathcal{C}\rightarrow A$ (called \textbf{counity}%
), such that%
\begin{equation*}
\dsum \varepsilon _{\mathcal{C}}(c_{1})c_{2}=c=\dsum c_{1}\varepsilon _{%
\mathcal{C}}(c_{2})\text{ for all }c\in \mathcal{C}.
\end{equation*}
\end{punto}

\begin{definition}
Let $(\mathcal{C},\Delta _{\mathcal{C}},\varepsilon _{\mathcal{C}})$ be an $%
A $-coring. An $A$-subbimodule $K\subseteq \mathcal{C}$ is a $\mathcal{C}$-%
\textbf{coideal}, iff%
\begin{equation*}
\Delta _{\mathcal{C}}(K)\subseteq K\otimes _{A}\mathcal{C}+\mathcal{C}%
\otimes _{A}K\text{ and }\varepsilon _{\mathcal{C}}(K)=0.
\end{equation*}
\end{definition}

\begin{punto}
A \emph{morphism of corings}%
\begin{equation*}
(\theta :\gamma ):(\mathcal{C},\Delta _{\mathcal{C}}:A)\rightarrow (\mathcal{%
D},\Delta _{\mathcal{D}}:B)
\end{equation*}%
consists of a morphism of unital $R$-algebras $\gamma :A\rightarrow B$ and
an $(A,A)$-bilinear map $\theta :\mathcal{C}\rightarrow \mathcal{D},$ such
that%
\begin{equation}
\chi _{(\mathcal{D},\mathcal{D})}^{(A,B)}\circ (\theta \otimes _{A}\theta
)\circ \Delta _{\mathcal{C}}=\Delta _{\mathcal{D}}\circ \theta .
\label{coring-mor}
\end{equation}%
In case the corings $\mathcal{C}$ and $\mathcal{D}$ are (left counital
(resp. right counital, counital), then we say $(\theta :\gamma )$ is a
morphism of \emph{left counital corings} (resp. \emph{right counital
corings, counital corings}), iff $\varepsilon _{\mathcal{D}}^{l}\circ \theta
=\gamma \circ \varepsilon _{\mathcal{C}}^{l}$ (resp. $\varepsilon _{\mathcal{%
D}}^{r}\circ \theta =\gamma \circ \varepsilon _{\mathcal{C}}^{r},$ $%
\varepsilon _{\mathcal{D}}\circ \theta =\gamma \circ \varepsilon _{\mathcal{C%
}}$).
\end{punto}

\begin{punto}
Let $(\mathcal{C},\Delta _{\mathcal{C}})$ be a (not necessarily counital) $A$%
-coring. Then

$(^{\ast }\mathcal{C},\ast _{l})$ is a (not necessarily unital) $A$-ring
with multiplication%
\begin{equation*}
(f\ast _{l}g)(c)=\dsum g(c_{1}f(c_{2}))\text{ for all }c\in \mathcal{C}\text{
and }f,g\in \text{ }^{\ast }\mathcal{C};
\end{equation*}

$(\mathcal{C}^{\ast },\ast _{r})$ is a (not necessarily unital) $A$-ring
with multiplication%
\begin{equation*}
(f\ast _{r}g)(c)=\dsum f(g(c_{1})c_{2})\text{ for all }c\in \mathcal{C}\text{
and }f,g\in \mathcal{C}^{\ast };
\end{equation*}

$(^{\ast }\mathcal{C}^{\ast },\ast )$ is a (not necessarily unital) $A$-ring
with multiplication%
\begin{equation*}
(f\ast g)(c)=\dsum g(c_{1})f(c_{2})\text{ for all }c\in \mathcal{C}\text{
and }f,g\in \text{ }^{\ast }\mathcal{C}^{\ast }.
\end{equation*}%
Notice that these multiplications are opposite to the ones used in \cite%
{BW-2003}.
\end{punto}

\begin{remark}
Let $(\mathcal{C},\Delta _{\mathcal{C}})$ be an $A$-coring.

if $(\mathcal{C},\Delta _{\mathcal{C}},\varepsilon _{\mathcal{C}}^{r})$ is a
right counital $A$-coring, then $(^{\ast }\mathcal{C},\ast _{l},\varepsilon
_{\mathcal{C}}^{r})$ is a left unital $A$-ring with left unity $\varepsilon
_{\mathcal{C}}^{r};$ the converse holds, if $_{A}\mathcal{C}$ is $A$%
-cogenerated.

if $(\mathcal{C},\Delta _{\mathcal{C}},\varepsilon _{\mathcal{C}}^{l})$ is a
left counital $A$-coring, then $(\mathcal{C}^{\ast },\ast _{r},\varepsilon _{%
\mathcal{C}}^{l})$ is a right unital $A$-ring with right unity $\varepsilon
_{\mathcal{C}}^{l};$ the converse holds, if $\mathcal{C}_{A}$ is $A$%
-cogenerated.

if $(\mathcal{C},\Delta _{\mathcal{C}},\varepsilon _{\mathcal{C}})$ is an $A$%
-coring, then $(^{\ast }\mathcal{C},\ast _{l},\varepsilon _{\mathcal{C}})$
(resp. $(\mathcal{C}^{\ast },\ast _{r},\varepsilon _{\mathcal{C}})$ and $%
(^{\ast }\mathcal{C}^{\ast },\ast ,\varepsilon _{\mathcal{C}})$) are $A$%
-rings; the converse holds if $_{A}\mathcal{C}$ (resp. $\mathcal{C}_{A},$ $%
_{A}\mathcal{C}_{A}$) is $A$-cogenerated.
\end{remark}

\begin{punto}
We call an $R$-coring with $rc=cr$ for all $c\in C$ and $r\in R$ a (non
necessarily counital) $R$\textbf{-coalgebra}. For any $R$-coalgebra $C,$ the
dual $R$-module $C^{\ast }:=\mathrm{Hom}_{R}(C,R)$ is an $R$-algebra with
multiplication given by the convolution product%
\begin{equation*}
(f\ast g)(c)=\dsum f(c_{1})g(c_{2})\text{ for all }f,g\in C^{\ast }\text{
and }c\in C.
\end{equation*}%
If $C$ is counital with counity $\varepsilon ,$ then the $R$-algebra $%
C^{\ast }$ is unital with unity $\varepsilon .$
\end{punto}

\begin{punto}
Let $(\mathcal{C},\Delta _{\mathcal{C}})$ be an $A$-coring. With a \textbf{%
right }$\mathcal{C}$\textbf{-comodule,} we mean a right $A$-module $M$ with
a right $A$-linear map $\varrho _{M}^{\mathcal{C}}:M\rightarrow M\otimes _{A}%
\mathcal{C},$ $m\mapsto \sum m_{<0>}\otimes _{A}m_{<1>}$ such that%
\begin{equation*}
(\mathrm{id}_{M}\otimes _{A}\Delta _{\mathcal{C}})\circ \varrho _{M}^{%
\mathcal{C}}=(\varrho _{M}^{\mathcal{C}}\otimes _{A}\mathrm{id}_{\mathcal{C}%
})\circ \varrho _{M}^{\mathcal{C}}.
\end{equation*}%
If $(\mathcal{C},\Delta _{\mathcal{C}},\varepsilon _{\mathcal{C}}^{r})$ is a
right counital $A$-coring we say $M$ is \textbf{counital}, iff%
\begin{equation*}
(\mathrm{id}_{M}\otimes _{A}\varepsilon _{\mathcal{C}}^{r})\circ \varrho
_{M}^{\mathcal{C}}\equiv \vartheta _{M}^{r},\text{ i.e. }\dsum
m_{<0>}\varepsilon _{\mathcal{C}}^{r}(m_{<1>})=m\text{ for all }m\in M.
\end{equation*}%
If $M$ and $N$ are right $\mathcal{C}$-comodules we say an $A$-linear map $%
f:M\rightarrow N$ is a \textbf{morphism of right }$\mathcal{C}$\textbf{%
-comodules} (called also \textbf{right }$\mathcal{C}$\textbf{-colinear})%
\textbf{, }iff $\varrho _{N}^{\mathcal{C}}\circ f=(f\otimes _{A}\mathrm{id}_{%
\mathcal{C}})\circ \varrho _{M}^{\mathcal{C}}.$ The category of right $%
\mathcal{C}$-comodules with right $\mathcal{C}$-colinear maps is denoted by $%
\mathbb{M}^{\mathcal{C}}.$ The category $^{\mathcal{C}}\mathbb{M}$ of \emph{%
left }$\mathcal{C}$\emph{-comodules} is defined analogously.
\end{punto}

\begin{punto}
Let $\mathcal{C}$ be an $A$-coring and $\mathcal{D}$ a $B$-coring. With a $(%
\mathcal{D},\mathcal{C})$-bicomodule, we mean a $(B,A)$-bimodule $M$ where $%
(M,\varrho _{M}^{\mathcal{D}})$ is a left $\mathcal{D}$-comodule and $%
(M,\varrho _{M}^{\mathcal{C}})$ is a right $\mathcal{C}$-comodule, such that
$\varrho _{M}^{\mathcal{D}}:M\rightarrow \mathcal{D}\otimes _{B}M$ is $%
\mathcal{C}$-colinear (equivalently, iff $\varrho _{M}^{\mathcal{C}%
}:M\rightarrow M\otimes _{A}\mathcal{C}$ is $\mathcal{D}$-colinear), i.e.%
\begin{equation*}
(\varrho _{M}^{\mathcal{D}}\otimes _{A}\mathrm{id}_{\mathcal{C}})\circ
\varrho _{M}^{\mathcal{C}}=(\mathrm{id}_{\mathcal{D}}\otimes _{B}\varrho
_{M}^{\mathcal{C}})\circ \varrho _{M}^{\mathcal{D}}.
\end{equation*}%
If $(M,\varrho _{M}^{\mathcal{D}},\varrho _{M}^{\mathcal{C}})$ and $%
(N,\varrho _{N}^{\mathcal{D}},\varrho _{N}^{\mathcal{C}})$ are $(\mathcal{D},%
\mathcal{C})$-bicomodules, then a $(B,A)$-bilinear map $f:M\rightarrow N$ is
called a \textbf{morphism of }$(\mathcal{D},\mathcal{C})$\textbf{-bicomodules%
} (or $(\mathcal{D},\mathcal{C})$\textbf{-bicolinear}), iff $f$ is left $%
\mathcal{D}$-colinear and right $\mathcal{C}$-colinear; and we denote the
set of such morphisms by $^{\mathcal{D}}\mathrm{Hom}^{\mathcal{C}}(M,N)$).
If $\mathcal{D}$ is a left counital $B$-coring with left counit $\varepsilon
_{\mathcal{D}}^{l}$ and $\mathcal{C}$ is a right counital $A$-coring with
right counit $\varepsilon _{\mathcal{C}}^{r}$ we say a $(\mathcal{D},%
\mathcal{C})$-bicomodule $M$ is \textbf{counital}, iff $M$ is counital as a
left $\mathcal{D}$-comodule and as a right $\mathcal{C}$-comodule, i.e. iff%
\begin{equation*}
\dsum \varepsilon _{\mathcal{D}}^{l}(m_{<-1>})m_{<0>}=m=\dsum
m_{<0>}\varepsilon _{\mathcal{C}}^{r}(m_{<1>})\text{ for all }m\in M.
\end{equation*}%
The category of $(\mathcal{D},\mathcal{C})$-bicomodules and $(\mathcal{D},%
\mathcal{C})$-bicolinear maps is denoted by $^{\mathcal{D}}\mathbb{M}^{%
\mathcal{C}}.$
\end{punto}

\begin{punto}
Let $\mathcal{C}=(\mathcal{C},\Delta _{\mathcal{C}}:A)$ be an $A$-coring.
Then $\mathcal{C}^{cop}=(\mathcal{C}^{cop},\Delta _{\mathcal{C}}^{\mathrm{tw}%
}:A^{op})$ is an $A^{op}$-coring, where%
\begin{equation*}
\Delta _{\mathcal{C}}^{\mathrm{tw}}:=\tau \circ \Delta _{\mathcal{C}}:%
\mathcal{C}\rightarrow \mathcal{C}^{cop}\otimes _{A^{op}}\mathcal{C}^{cop},%
\text{ }c\mapsto \dsum c_{2}\otimes _{A^{op}}c_{1}.
\end{equation*}%
Moreover, the category or right (left) $\mathcal{C}$-comodules is isomorphic
to the category of left (right) $\mathcal{C}$-comodules; and the category of
$(\mathcal{D},\mathcal{C})$-bicomodules, where $\mathcal{D}$ is a $B$%
-coring, is isomorphic to the category of $(\mathcal{C}^{cop},\mathcal{D}%
^{cop})$-bicomodules.
\end{punto}

\begin{remark}
If $\mathcal{C}$ is left counital (resp. right counital, counital) $A$%
-coring, then the $A^{op}$-coring $\mathcal{C}^{cop}$ is right counital
(resp. left counital, counital).
\end{remark}

\qquad A slight modification of the proofs of \cite{Abu-2003} (or \cite%
{BW-2003}) yields the following

\begin{proposition}
\label{sg} Let $\mathcal{C}\ $be an $A$-coring.

\begin{enumerate}
\item Assume $\mathcal{C}$ is a right counital $A$-coring. Then $_{A}%
\mathcal{C}$ is locally projective $\Leftrightarrow $ $\mathbb{M}^{\mathcal{C%
}}\simeq \sigma \lbrack \mathcal{C}_{^{\ast }\mathcal{C}}]\Leftrightarrow $ $%
\mathbb{M}^{\mathcal{C}}\subseteq \mathbb{M}_{^{\ast }\mathcal{C}}$ is a
full subcategory;

\item Assume $\mathcal{C}$ is a left counital $A$-coring. Then $\mathcal{C}%
_{A}$ is locally projective $\Leftrightarrow $ $^{\mathcal{C}}\mathbb{M}%
\simeq \sigma \lbrack _{\mathcal{C}^{\ast }}\mathcal{C}]\Leftrightarrow $ $^{%
\mathcal{C}}\mathbb{M}\subseteq $ $_{\mathcal{C}^{\ast }}\mathbb{M}$ is a
full subcategory.
\end{enumerate}
\end{proposition}

\subsection*{Dorroh Corings}

\qquad Inspired by the notion of Dorroh rings, obtained by associating a
unity to a non-unital ring, Vercruysse introduced in \cite{Ver-2006} the
notion of Dorroh corings obtained by associating a counity to a non-counital
coring:

\begin{theorem}
\label{Dorr}\emph{(\cite[Theorem 6.1.]{Ver-2006}) }Let $\mathcal{C}$ be an $%
A $-coring. Then there exists a counital $A$-coring $\widehat{\mathcal{C}},$
such that:

\begin{enumerate}
\item There is a surjective morphism of $A$-corings $\pi :\widehat{\mathcal{C%
}}\rightarrow \mathcal{C};$

\item $\mathcal{C}$ is isomorphic to a coideal of $\widehat{\mathcal{C}};$

\item There is an injective morphism of $A$-corings $\iota :A\rightarrow
\widehat{\mathcal{C}};$

\item The category of \emph{not necessarily counital }right \emph{(}left%
\emph{) }$\mathcal{C}$-comodules is isomorphism to the category of counital
right \emph{(}left\emph{)} $\widehat{\mathcal{C}}$-comodules.
\end{enumerate}
\end{theorem}

\begin{definition}
Let $\mathcal{C}$ be an $A$-coring. The $A$-coring $\widehat{\mathcal{C}}$
satisfying the equivalent conditions of Theorem \ref{Dorr} is called the
\textbf{Dorroh coring} associated to $\mathcal{C}.$
\end{definition}

\begin{punto}
\label{hat-str}Let $A$ be a ring and $\mathcal{C}$ an $A$-coring. The
associated Dorroh $A$-coring $\widehat{\mathcal{C}}$ is constructed as
follows: $\widehat{\mathcal{C}}:=\mathcal{C}\times A$ with bimodule
structure given by
\begin{equation*}
a^{\prime }(c,a)a^{\prime \prime }:=(a^{\prime }ca^{\prime \prime
},a^{\prime }aa^{\prime \prime })\text{ for all }a,a^{\prime },a^{\prime
\prime }\in A\text{ and }c\in \mathcal{C}.
\end{equation*}%
The comultiplication and the counit of $\widehat{\mathcal{C}}$ are%
\begin{equation*}
\begin{tabular}{lll}
$\Delta _{\widehat{\mathcal{C}}}(c,a)$ & $:=$ & $\dsum (c_{1},0)\otimes
_{A}(c_{2},0)+(0,1_{A})\otimes _{A}(c,a)+(c,a)\otimes
_{A}(0,1_{A})-(0,a)\otimes _{A}(0,1_{A});$ \\
$\varepsilon _{\widehat{\mathcal{C}}}(c,a)$ & $:=$ & $a.$%
\end{tabular}%
\end{equation*}%
If $(M,\varrho _{M}^{\mathcal{C}})$ is a right $\mathcal{C}$-comodule, then $%
(M,\varrho _{M}^{\widehat{\mathcal{C}}})$ is a \emph{counital} right $%
\widehat{\mathcal{C}}$-comodule, where%
\begin{equation*}
\varrho _{M}^{\widehat{\mathcal{C}}}:M\rightarrow M\otimes _{A}\widehat{%
\mathcal{C}},\text{ }m\mapsto \dsum m_{<0>}\otimes _{A}(m_{<1>},0)+m\otimes
_{A}(0,1_{A}).
\end{equation*}%
On the other hand, if $(N,\varrho _{N}^{\mathcal{C}})$ is a left $\mathcal{C}
$-comodule, then $(N,\varrho _{N}^{\widehat{\mathcal{C}}})$ is a \emph{%
counital} $\widehat{\mathcal{C}}$-comodule, where%
\begin{equation*}
\varrho _{N}^{\widehat{\mathcal{C}}}:N\rightarrow \widehat{\mathcal{C}}%
\otimes _{A}N,\text{ }n\mapsto \dsum (n_{<-1>},0)\otimes
_{A}n_{<0>}+(0,1_{A})\otimes _{A}n.
\end{equation*}
\end{punto}

\qquad Similar to Theorem \ref{Dorr}, we get the following result:

\begin{theorem}
\label{Bi-Dorr}Let $\mathcal{D}$ and $\mathcal{C}$ be corings, and consider
the associated counital Dorroh corings $\widehat{\mathcal{D}}$ and $\widehat{%
\mathcal{C}}.$ Then the category of \emph{(}not necessarily counital\emph{) }%
$(\mathcal{D},\mathcal{C})$-bicomodules is isomorphic to the category of
counital $(\widehat{\mathcal{D}},\widehat{\mathcal{C}})$-bicomodules.
\end{theorem}

\section{Every Coalgebra is a Coring}

\bigskip

\qquad In this section we show that every $R$-coalgebra $C$ is a (possibly
non-counital) $C^{\bullet }$-coring, where $C^{\bullet }:=C^{\ast op}$ is
the opposite of the dual $R$-algebra.

\bigskip

\begin{punto}
Let $(C,\Delta _{C},\varepsilon _{C})$ be an $R$-coalgebra and consider the
\emph{opposite dual }$R$\emph{-algebra }$C^{\bullet }:=((C^{\ast
})^{op},\bullet ,\varepsilon _{C}),$ where%
\begin{equation*}
(f\bullet g)(c)=\dsum g(c_{1})f(c_{2})\text{ for all }f,g\in C^{\bullet }%
\text{ and }c\in C.
\end{equation*}%
Then $C$ is clearly a $(C^{\bullet },C^{\bullet })$-bimodule, with left and
right $C^{\bullet }$-actions given by%
\begin{equation*}
f\rightharpoondown c:=\dsum f(c_{1})c_{2}\text{ and }c\leftharpoondown
g:=\dsum c_{1}g(c_{2})\text{ for all }f,g\in C^{\bullet }\text{ and }c\in C.
\end{equation*}
\end{punto}

\begin{proposition}
Let $(C,\Delta _{C},\varepsilon _{C})$ be an $R$-coalgebra and consider the $%
R$-linear morphism%
\begin{equation*}
\eta _{C^{\bullet }}:R\rightarrow C^{\bullet },\text{ }r\mapsto \lbrack
c\mapsto r\varepsilon _{C}(c)].
\end{equation*}

\begin{enumerate}
\item $\mathbf{C}:=(\mathbf{C}:C^{\bullet })$ is a \emph{(}not necessarily
counital\emph{)} coring, where $\mathbf{C}=C$ is canonically a $(C^{\bullet
},C^{\bullet })$-bimodule, and with comultiplication%
\begin{equation*}
\Delta _{\mathbf{C}}:=\chi _{(C,C)}^{(R,C^{\bullet })}\circ \Delta _{C}:%
\mathbf{C}\rightarrow \mathbf{C}\otimes _{C^{\bullet }}\mathbf{C},\text{ }%
c\mapsto \dsum c_{1}\otimes _{C^{\bullet }}c_{2}\text{ for every }c\in C.
\end{equation*}%
Moreover, we have a morphism of corings%
\begin{equation*}
(\mathrm{id}:\eta _{C^{\bullet }}):(C:R)\rightarrow (\mathbf{C}:C^{\bullet
}).
\end{equation*}

\item $\widehat{\mathbf{C}}:=\mathbf{C}\times C^{\bullet }$ is a counital $%
C^{\bullet }$-coring with canonical $(C^{\bullet },C^{\bullet })$-bimodule
structure, and with comultiplication and counit%
\begin{equation*}
\begin{tabular}{lll}
$\Delta _{\widehat{\mathbf{C}}}(c,f)$ & $:=$ & $\dsum (c_{1},0)\otimes
_{C^{\bullet }}(c_{2},0)+(0,\varepsilon _{C})\otimes _{C^{\bullet
}}(c,f)+(c,f)\otimes _{C^{\bullet }}(0,\varepsilon _{C})$ \\
&  & $-(0,f)\otimes _{C^{\bullet }}(0,\varepsilon _{C});$ \\
$\varepsilon _{\widehat{\mathbf{C}}}(c,f)$ & $:=$ & $f.$%
\end{tabular}%
\end{equation*}

\item The category of \emph{not necessarily counital }right \emph{(}left%
\emph{)} $\mathbf{C}$-comodules is isomorphic to the category of \emph{%
counital} right \emph{(}left\emph{)} $\widehat{\mathbf{C}}$-comodules.

\item If $D$ is any $R$-coalgebra and $(\mathbf{D}:D^{\bullet })$ is the
corresponding $D^{\ast op}$-coring, then the category of \emph{not
necessarily counital }$(\mathbf{D},\mathbf{C})$-bicomodules is isomorphic to
the category of \emph{counital} $(\widehat{\mathbf{D}},\widehat{\mathbf{C}})$%
-bicomodules.
\end{enumerate}
\end{proposition}

\begin{Beweis}
Since $\Delta _{C}$ is coassociative, it is clear that $\Delta _{\mathbf{C}%
}:=\chi _{(C,C)}^{(R,C^{\bullet })}\circ \Delta _{C}$ is coassociative as
well. Moreover, we have for all $c\in C$ and $f,g\in C^{\bullet }:$%
\begin{equation*}
\Delta _{\mathbf{C}}(f\rightharpoondown c)=\Delta _{\mathbf{C}}(\dsum
f(c_{1})c_{2})=\dsum f(c_{1})c_{2}\otimes _{C^{\bullet
}}c_{3}=f\rightharpoondown (\dsum c_{1}\otimes _{C^{\bullet
}}c_{2})=f\rightharpoondown \Delta _{\mathbf{C}}(c)
\end{equation*}%
and%
\begin{equation*}
\Delta _{\mathbf{C}}(c\leftharpoondown g)=\Delta _{\mathbf{C}}(\dsum
c_{1}g(c_{2}))=\dsum c_{1}\otimes _{C^{\bullet }}c_{2}g(c_{3})=(\dsum
c_{1}\otimes _{C^{\bullet }}c_{2})\leftharpoondown g=\Delta _{\mathbf{C}%
}(c)\leftharpoondown g,
\end{equation*}%
whence $\Delta _{\mathbf{C}}$ is $(C^{\bullet },C^{\bullet })$-bilinear. It
follows then that $(\mathbf{C},\Delta _{\mathbf{C}})$ is a (not necessarily
counital) $A$-coring. The remaining results follow from Theorems \ref{Dorr}, %
\ref{Bi-Dorr} and \ref{hat-str}.$\blacksquare $
\end{Beweis}

\begin{remarks}
(\cite[8.10.]{BW-2003}) Let $(C,\Delta _{C},\varepsilon _{C})$ be an $R$%
-coalgebra and assume $_{R}C$ is locally projective (so that $\mathbb{M}%
^{C}\simeq \sigma \lbrack _{C^{\ast }}C]=\sigma \lbrack C_{C^{\bullet }}]$
and $^{C}\mathbb{M}\simeq \sigma \lbrack C_{C^{\ast }}]=\sigma \lbrack
_{C^{\bullet }}C]$).

\begin{enumerate}
\item If $_{C^{\bullet }}C$ ($C_{C^{\bullet }}$) is locally projective, then
$C$ is a generator in $\mathbb{M}^{C}$ (in $^{C}\mathbb{M}$).

\item If $_{C^{\bullet }}C$ and $C_{C^{\bullet }}$ are locally projective,
then the functors $\mathrm{Rat}^{C}(-):\mathbb{M}_{C^{\bullet }}\rightarrow
\mathbb{M}^{C}$ and $^{C}\mathrm{Rat}(-):$ $_{C^{\bullet }}\mathbb{M}%
\rightarrow $ $^{C}\mathbb{M}$ are exact.
\end{enumerate}
\end{remarks}

\begin{corollary}
Let $(C,\Delta _{C},\varepsilon _{C})$ be an $R$-coalgebra and consider the
morphism of $R$-algebras%
\begin{equation*}
\eta _{C^{\ast }}:R\rightarrow C^{\ast },\text{ }r\mapsto \lbrack c\mapsto
r\varepsilon _{C}(c)].
\end{equation*}

\begin{enumerate}
\item $\mathbf{C^{cop}}:=(C:C^{\ast })$ is a $C^{\ast }$-coring with the
canonical $(C^{\ast },C^{\ast })$-bimodule structure, and with
comultiplication given by%
\begin{equation*}
\Delta _{\mathbf{C^{cop}}}:=\chi _{(C,C)}^{(R,C^{\ast })}\circ \Delta ^{%
\mathrm{tw}}:\mathbf{C^{cop}}\rightarrow \mathbf{C^{cop}}\otimes _{C^{\ast }}%
\mathbf{C^{cop}},\text{ }c\mapsto \dsum c_{2}\otimes _{C^{\ast }}c_{1}\text{
for every }c\in C.
\end{equation*}%
Moreover, we have a morphism of corings%
\begin{equation*}
(\mathrm{id}:\eta _{C^{\ast }}):(C:R)\rightarrow (\mathbf{C^{cop}}:C^{\ast
}).
\end{equation*}

\item $\widehat{\mathbf{C^{cop}}}:=\mathbf{C^{cop}}\times C^{\ast }$ is a
counital $C^{\ast }$-coring with canonical $(C^{\ast },C^{\ast })$-bimodule
structure, and with comultiplication and counit given by%
\begin{equation*}
\begin{tabular}{lll}
$\mathbf{\Delta }_{\widehat{\mathbf{C^{cop}}}}(c,f)$ & $:=$ & $\dsum
(c_{2},0)\otimes _{C^{\ast }}(c_{1},0)+(0,\varepsilon _{C})\otimes _{C^{\ast
}}(c,f)+(c,f)\otimes _{C^{\ast }}(0,\varepsilon _{C})$ \\
&  & $-(0,f)\otimes _{C^{\ast }}(0,\varepsilon _{C});$ \\
$\varepsilon _{\widehat{\mathbf{C^{cop}}}}(c,f)$ & $:=$ & $f.$%
\end{tabular}%
\end{equation*}

\item The category of \emph{not necessarily counital }right (left) $\mathbf{%
C^{cop}}$-comodules is isomorphic to the category of \emph{counital} right
(left) $\widehat{\mathbf{C^{cop}}}$-comodules; and the category of \emph{not
necessarily counital} $(\mathbf{C^{cop}},\mathbf{C^{cop}})$-bicomodules is
isomorphic to the category of \emph{counital} $(\widehat{\mathbf{C^{cop}}},%
\widehat{\mathbf{C^{cop}}})$-bicomodules.
\end{enumerate}
\end{corollary}

\begin{remark}
Let $(C,\Delta _{C},\varepsilon _{C})$ be an $R$-coalgebra and consider the $%
R$-linear map%
\begin{equation*}
\overline{\varepsilon }:\mathbf{C}\rightarrow C^{\bullet },\text{ }c\mapsto
\varepsilon _{C}(c)\text{ }\varepsilon _{C}(-).
\end{equation*}%
For any $f\in C^{\bullet }$ and $c,d\in C$ we have%
\begin{equation*}
\overline{\varepsilon }(f\rightharpoondown c)(d)=\overline{\varepsilon }%
(\dsum f(c_{1})c_{2})(d)=\dsum f(c_{1})\varepsilon _{C}(c_{2})\varepsilon
_{C}(d)=f(\dsum c_{1}\varepsilon _{C}(c_{2}))\varepsilon
_{C}(d)=f(c)\varepsilon _{C}(d),
\end{equation*}%
while%
\begin{equation*}
\lbrack f\bullet \overline{\varepsilon }(c)](d)=\dsum [\overline{\varepsilon
}(c)(d_{1})]f(d_{2})=\dsum \varepsilon _{C}(c)\varepsilon
_{C}(d_{1})f(d_{2})=\varepsilon _{C}(c)f(\dsum \varepsilon
_{C}(d_{1})d_{2})=\varepsilon _{C}(c)f(d),
\end{equation*}%
whence $\overline{\varepsilon }$ is not left $C^{\bullet }$-linear.
Similarly, one can show that $\overline{\varepsilon }$ is not right $%
C^{\bullet }$-linear. Hence $\varepsilon $ cannot serve as a left (right)
counit for the $C^{\bullet }$-coring $\mathbf{C},$ or for the $C^{\ast }$%
-coring $\mathbf{C^{cop}.}$ We conclude then that these corings are \emph{%
not necessarily} counital, and so the Dorroh corings above are non-trivial.
\end{remark}

\begin{definition}
A (not necessarily unital) $R$-algebra $(A,\mu _{\mathbf{A}})$ is said to be
\textbf{separable}, iff $\mu _{\mathbf{A}}:\mathbf{A}\otimes _{R}\mathbf{A}%
\rightarrow \mathbf{A}$ has a \emph{section} (i.e. there exists a $(\mathbf{A%
},\mathbf{A})$-bilinear map $\delta _{\mathbf{A}}:\mathbf{A}\rightarrow
\mathbf{A}\otimes _{R}\mathbf{A},$ such that $\mu _{\mathbf{A}}\circ \delta
_{\mathbf{A}}=\mathrm{id}_{\mathbf{A}}.$
\end{definition}

\begin{definition}
A \emph{separable} \emph{measuring }$R$\emph{-pairing} $P=(A,C;\kappa _{P})$
consists of a \emph{separable} (not necessarily unital) $R$-algebra $\mathbf{%
A}$ and a (not necessarily counital) \emph{coseparable }$R$-coalgebra $C$
along with a morphism of $R$-algebras $\kappa _{P}:\mathbf{A}\rightarrow
C^{\ast }.$
\end{definition}

\begin{proposition}
\label{C-separable}Let $(C,\Delta _{C})$ be a \emph{(not necessarily
counital)} $R$-coalgebra. If the associated $C^{\bullet }$-coring $\mathbf{C}%
=(C:C^{\bullet })\mathbf{\ }$is left \emph{(}right\emph{)} counital, then

\begin{enumerate}
\item $C$ is a separable \emph{(not necessarily unital) }$R$-algebra with
multiplication%
\begin{equation*}
\mu _{C}:C\otimes _{R}C\rightarrow C,\text{ }c\otimes _{R}c^{\prime }\mapsto
\varepsilon ^{l}(c)\rightharpoondown c^{\prime }.
\end{equation*}

\item $C$ is a coseparable $R$-coalgebra.

\item $P=(C^{op},C)$ is a separable measuring $R$-pairing \emph{(}where $%
C^{op}$ is the opposite $R$-algebra\emph{)}.
\end{enumerate}
\end{proposition}

\begin{Beweis}
Let $(C,\Delta _{C})$ be a (\emph{not necessarily counital}) $R$-coalgebra
and assume the associated $C^{\bullet }$-coring $\mathbf{C}=(C:C^{\bullet })$
to be have a left counit $\varepsilon _{\mathbf{C}}^{l}:C\rightarrow
C^{\bullet }$ (if $\mathbf{C}$ has a right counit $\varepsilon _{\mathbf{C}%
}^{r},$ then the results can be proved analogously).

\begin{enumerate}
\item We have for all $c,d,e\in C:$%
\begin{equation*}
\begin{tabular}{lll}
$(\mu _{C}\circ (\mu _{C}\otimes _{R}\mathrm{id}_{C}))(c\otimes _{R}d\otimes
_{R}e)$ & $=$ & $\mu _{C}(\varepsilon _{\mathbf{C}}^{l}(c)\rightharpoondown
d\otimes _{R}e)$ \\
& $=$ & $\varepsilon ^{l}(\varepsilon _{\mathbf{C}}^{l}(c)\rightharpoondown
d)\rightharpoondown e$ \\
& $=$ & $(\varepsilon _{\mathbf{C}}^{l}(c)\cdot \varepsilon _{\mathbf{C}%
}^{l}(d))\rightharpoondown e$ \\
& $=$ & $\varepsilon _{\mathbf{C}}^{l}(c)\rightharpoondown (\varepsilon _{%
\mathbf{C}}^{l}(d)\rightharpoondown e)$ \\
& $=$ & $\mu _{C}(c\otimes _{R}\varepsilon _{\mathbf{C}}^{l}(d)%
\rightharpoondown e)$ \\
& $=$ & $(\mu _{C}\circ (\mathrm{id}_{C}\otimes _{R}\mu _{C}))(c\otimes
_{R}d\otimes _{R}e),$%
\end{tabular}%
\end{equation*}%
whence $(C,\mu _{C})$ is an $R$-algebra. Notice also that for all $%
c,c^{\prime }\in C$ we have%
\begin{equation*}
\begin{tabular}{lllll}
$\Delta _{C}(cc^{\prime })$ & $=$ & $\Delta _{C}(\varepsilon _{\mathbf{C}%
}^{l}(c)\rightharpoondown c^{\prime })$ & $=$ & $\varepsilon _{\mathbf{C}%
}^{l}(c)\rightharpoondown \Delta _{C}(c^{\prime })$ \\
& $=$ & $\varepsilon _{\mathbf{C}}^{l}(c)\rightharpoondown \dsum
c_{1}^{\prime }\otimes _{R}c_{2}^{\prime }$ & $=$ & $\dsum \varepsilon _{%
\mathbf{C}}^{l}(c)\rightharpoondown c_{1}^{\prime }\otimes _{R}c_{2}^{\prime
}$ \\
& $=$ & $\dsum cc_{1}^{\prime }\otimes _{R}c_{2}^{\prime }$ & $=$ & $c\Delta
_{C}(c^{\prime }),$%
\end{tabular}%
\end{equation*}%
i.e. $\Delta _{C}$ is left $C$-linear. Similarly, one can show that $\Delta
_{C}$ is right $C$-linear. Moreover, for all $c\in C$ we have%
\begin{equation*}
(\mu _{C}\circ \Delta _{C})(c)=\mu _{C}(\dsum c_{1}\otimes _{R}c_{2})=\dsum
\varepsilon _{\mathbf{C}}^{l}(c_{1})\rightharpoondown c_{2}=c,
\end{equation*}%
i.e. $\Delta _{C}:C\rightarrow C\otimes _{R}C$ is a section of $\mu
_{C}:C\otimes _{R}C\rightarrow C.$

\item This follows by \cite[26.8.]{BW-2003}.

\item Consider%
\begin{equation*}
\kappa _{P}:C^{op}\rightarrow C^{\ast },\text{ }c\mapsto \lbrack d\mapsto
\varepsilon _{\mathbf{C}}^{l}(c)(d)].
\end{equation*}%
For all $c,c^{\prime },d\in C$ we have%
\begin{equation*}
\kappa _{P}(c\bullet ^{op}c^{\prime })=\varepsilon _{\mathbf{C}%
}^{l}(c^{\prime }c)=\varepsilon _{\mathbf{C}}^{l}(\varepsilon _{\mathbf{C}%
}^{l}(c^{\prime })\rightharpoondown c)=\varepsilon _{\mathbf{C}%
}^{l}(c^{\prime })\bullet \varepsilon ^{l}(c)=\kappa _{P}(c^{\prime
})\bullet \kappa _{P}(c)=\kappa _{P}(c)\ast \kappa _{P}(c^{\prime }),
\end{equation*}%
i.e. $\kappa _{P}$ is a morphism of $R$-algebras. Whence $P=(C^{op},C)$ is a
measuring $R$-pairing.$\blacksquare $
\end{enumerate}
\end{Beweis}

\section{Coseparable Corings}

\qquad Throughout this section, $C=(C,\Delta _{C},\varepsilon _{C})$ is a
\emph{counital} $R$-coalgebra and $C^{cop}:=(C,\Delta _{C}^{\mathrm{tw}%
},\varepsilon _{C})$ is its opposite $R$-coalgebra.

\begin{definition}
We say that a (not necessarily counital) $A$-coring $(\mathcal{C},\Delta _{%
\mathcal{C}})$ is \textbf{coseparable}, iff the structure map $\Delta _{%
\mathcal{C}}:\mathcal{C}\rightarrow \mathcal{C}\otimes _{A}\mathcal{C}$
splits as a $(\mathcal{C},\mathcal{C})$-bicomodule morphism, i.e. iff there
exists an $(A,A)$-bilinear map $\pi :\mathcal{C}\otimes _{A}\mathcal{C}%
\rightarrow \mathcal{C}$ such that%
\begin{equation*}
\pi \circ \Delta _{\mathcal{C}}=\mathrm{id}_{\mathcal{C}}\text{ and }(%
\mathrm{id}_{\mathcal{C}}\otimes _{A}\pi )\circ (\Delta _{\mathcal{C}%
}\otimes _{A}\mathrm{id}_{\mathcal{C}})=\Delta _{\mathcal{C}}\circ \pi =(\pi
\otimes _{A}id_{\mathcal{C}})\circ (\mathrm{id}_{\mathcal{C}}\otimes
_{A}\Delta _{\mathcal{C}}).
\end{equation*}
\end{definition}

\begin{definition}
Let $(\mathcal{C},\Delta _{\mathcal{C}},\varepsilon _{\mathcal{C}})$ be a
counital $A$-coring. A \textbf{cointegral in }$\mathcal{C}$ is an $(A,A)$%
-bilinear map $\gamma :\mathcal{C}\otimes _{A}\mathcal{C}\rightarrow A,$
such that for all $c,c^{\prime }\in \mathcal{C}$ we have%
\begin{equation*}
\dsum \gamma (c\otimes _{A}c_{1}^{\prime })c_{2}^{\prime }=\dsum c_{1}\gamma
(c_{2}\otimes _{A}c^{\prime })\text{ and }\dsum \gamma (c_{1}\otimes
_{A}c_{2})=\varepsilon _{\mathcal{C}}(c).
\end{equation*}
\end{definition}

\begin{lemma}
The mapping $\pi \mapsto \varepsilon _{\mathcal{C}}\circ \pi $ gives a 1-1
correspondence%
\begin{equation*}
\{\pi \in \text{ }^{\mathcal{C}}\mathrm{Hom}^{\mathcal{C}}(\mathcal{C}%
\otimes _{A}\mathcal{C},\mathcal{C})\mid \pi \circ \Delta _{\mathcal{C}}=%
\mathrm{id}_{\mathcal{C}}\}\longleftrightarrow \{\gamma \in \mathrm{Hom}%
_{(A,A)}(\mathcal{C}\otimes _{A}\mathcal{C},A)\mid \gamma \text{ cointegral
in }\mathcal{C}\},
\end{equation*}%
with inverse mapping $\gamma \mapsto \lbrack c\otimes _{A}c^{\prime }\mapsto
\dsum c_{1}\gamma (c_{2}\otimes _{A}c^{\prime })].$
\end{lemma}

\qquad \emph{Separable functors }were introduced first by N\u{a}st\u{a}sescu
et. al. \cite{NvdBvO1989} in the context of graded rings and were
investigated by several authors including Rafael \cite{Raf1990}. In \cite%
{C-IG-TN1997} (resp. \cite{G-T-2002}, \cite{G-TL2003}), G\'{o}%
mez-Torrecillas et. al. applied separable functors to study several functors
between categories of comodules for coalgebras over base fields (corings
over \emph{non-necessarily commutative} ground rings).

\begin{definition}
A covariant functor $\mathcal{F}:\mathcal{A}\rightarrow \mathcal{B}$ is said
to be \textbf{separable}, iff the functorial morphism%
\begin{equation*}
\mathcal{F}:\mathrm{Mor}_{\mathbf{A}}(-,-)\rightarrow \mathrm{Mor}_{\mathbf{B%
}}(F(-),F(-))
\end{equation*}%
is a \emph{functorial coretraction}, i.e. iff there exists a functorial
morphism%
\begin{equation*}
\mathcal{G}:\mathrm{Mor}_{\mathbf{B}}(F(-),F(-))\rightarrow \mathrm{Mor}_{%
\mathbf{A}}(-,-),
\end{equation*}%
such that $\mathcal{G}\circ \mathcal{F}=\mathrm{id}_{\mathrm{Mor}_{\mathbf{A}%
}(-,-)}.$
\end{definition}

\begin{lemma}
\label{C-Cop}The following are equivalent:

\begin{enumerate}
\item $C$ is coseparable;

\item there exists an $R$-linear map $\delta :C\otimes _{R}C\rightarrow R$
such that%
\begin{equation*}
\delta \circ \Delta _{C}=\varepsilon _{C}\text{ and }(\mathrm{id}_{C}\otimes
_{R}\delta )\circ (\Delta _{C}\otimes _{R}\mathrm{id}_{C})=(\delta \otimes
_{R}\mathrm{id}_{C})\circ \mathrm{(id}_{C}\otimes _{R}\Delta _{C}).
\end{equation*}

\item The forgetful functor $(-)_{R}:\mathbb{M}^{C}\rightarrow \mathbb{M}%
_{R} $ is separable;

\item The forgetful functor $_{R}(-)_{R}:$ $^{C}\mathbb{M}^{C}\rightarrow $ $%
_{R}\mathbb{M}_{R}$ is separable;

\item The forgetful functor $_{R}(-):$ $^{C}\mathbb{M}\rightarrow $ $_{R}%
\mathbb{M}$ is separable;

\item there exists an $R$-linear map $\widetilde{\delta }:C\otimes
_{R}C\rightarrow R$ such that%
\begin{equation*}
\widetilde{\delta }\circ \Delta _{C}^{\mathrm{tw}}=\varepsilon _{C}\text{
and }(\mathrm{id}_{C}\otimes _{R}\widetilde{\delta })\circ (\Delta _{C}^{%
\mathrm{tw}}\otimes _{R}id_{C})=(\widetilde{\delta }\otimes _{R}\mathrm{id}%
_{C})\circ (\mathrm{id}_{C}\otimes _{R}\Delta _{C}^{\mathrm{tw}});
\end{equation*}

\item $C^{cop}$ is coseparable.
\end{enumerate}
\end{lemma}

\begin{Beweis}
The statements $1$-$4$ are equivalent by \cite[3.29]{BW-2003}. The
conditions $4$-$7$ are equivalent, by restating the equivalent conditions $1$%
-$4$ for the opposite $R$-coalgebra $C^{cop}$ (notice that we have
isomorphisms of categories $^{C}\mathbb{M}^{C}\simeq $ $^{C^{cop}}\mathbb{M}%
^{C^{cop}}$ and $^{C}\mathbb{M}\simeq \mathbb{M}^{C^{op}}.\blacksquare $
\end{Beweis}

\qquad Next, we prove (under weaker assumptions) a technical lemma that is
similar to \cite[6.4.]{BW-2003}:

\begin{lemma}
\label{balanced}\ Let%
\begin{equation*}
\gamma :C\times C\rightarrow R,\text{ }(c,d)\mapsto <c,d>
\end{equation*}%
be an $R$-balanced form, and consider the associated $R$-linear maps:%
\begin{equation*}
\begin{tabular}{llllllll}
$\widetilde{\gamma }$ & $:$ & $C\otimes _{R}C$ & $\rightarrow $ & $R,$ & $%
c\otimes _{R}d$ & $\mapsto $ & $<c,d>;$ \\
$\gamma ^{l}$ & $:$ & $\mathbf{C}$ & $\rightarrow $ & $C^{\bullet },$ & $d$
& $\mapsto $ & $<-,d>;$ \\
$\gamma ^{r}$ & $:$ & $\mathbf{C}$ & $\rightarrow $ & $C^{\bullet },$ & $c$
& $\mapsto $ & $<c,->.$%
\end{tabular}%
\end{equation*}%
The following are equivalent:

\begin{enumerate}
\item $\gamma $ is $C^{\bullet }$-balanced;

\item $\gamma ^{l}:\mathbf{C}\rightarrow C^{\bullet }$ is left $C^{\bullet }$%
-linear;

\item $\gamma ^{r}:\mathbf{C}\rightarrow C^{\bullet }$ is right $C^{\bullet
} $-linear;

\item $\gamma $ factors through some $R$-linear map $\delta :C\otimes
_{C^{\bullet }}C\rightarrow R.$

If moreover $C$ is $R$-cogenerated, then the statements above are equivalent
to:

\item $(\mathrm{id}_{C}\otimes _{R}\widetilde{\gamma })\circ (\Delta _{C}^{%
\mathrm{tw}}\otimes _{R}\mathrm{id}_{C})=(\widetilde{\gamma }\otimes _{R}%
\mathrm{id}_{C})\circ (\mathrm{id}_{C}\otimes _{R}\Delta _{C}^{\mathrm{tw}%
}). $
\end{enumerate}
\end{lemma}

\begin{Beweis}
$(1)\Rightarrow (2):$ Assume that $\gamma $ is $C^{\bullet }$-balanced. For
any $f\in C^{\bullet }$ and $c,d\in \mathbf{C}$ we have%
\begin{equation*}
\begin{tabular}{lllll}
$\gamma ^{l}(f\rightharpoondown d)(c)$ & $=$ & $<c,f\rightharpoondown d>$ & $%
=$ & $<c\leftharpoondown f,d>$ \\
& $=$ & $\gamma ^{l}(d)(c\leftharpoondown f)$ & $=$ & $\gamma ^{l}(d)(\dsum
c_{1}f(c_{2}))$ \\
& $=$ & $(f\ast _{l}\gamma ^{l}(d))(c),$ &  &
\end{tabular}%
\end{equation*}%
i.e. $\gamma ^{l}$ is left $C^{\bullet }$-linear.

$(2)\Rightarrow (1):$ The result follows by rearranging the equalities above.

$(1)\Leftrightarrow (3):$ Follows by an argument similar to the one in
proving $(1)\Leftrightarrow (2).$

$(1)\Leftrightarrow (4):$ The equivalence follows by the definition of
tensor products.

$(1)\Rightarrow (5):$ Let $C$ be $R$-cogenerated and assume that $\gamma $
is $C^{\bullet }$-balanced. Then we have for every $c,d\in C$ and $f\in
C^{\ast }:$%
\begin{equation*}
\begin{tabular}{lll}
$f[(\mathrm{id}_{C}\otimes _{R}\widetilde{\gamma })\circ (\Delta _{C}^{%
\mathrm{tw}}\otimes _{R}\mathrm{id}_{C})(c\otimes _{R}d)]$ & $=$ & $\dsum
f(c_{2})<c_{1},d>$ \\
& $=$ & $<c\leftharpoondown f,d>$ \\
& $=$ & $<c,f\rightharpoondown d>$ \\
& $=$ & $\dsum <c,d_{2}>f(d_{1})$ \\
& $=$ & $f[(\widetilde{\gamma }\otimes _{R}\mathrm{id}_{C})\circ (\mathrm{id}%
_{C}\otimes _{R}\Delta _{C}^{\mathrm{tw}})(c\otimes _{R}d)].$%
\end{tabular}%
\end{equation*}%
Since $C$ is $R$-cogenerated, we obtain \textquotedblleft
5\textquotedblright .

$(5)\Rightarrow (1):$ Assuming \textquotedblleft 5\textquotedblright , we
get \textquotedblleft 1\textquotedblright\ by a similar argument to the one
above (notice that there is no need to assume $C$ is $R$-cogenerated to
prove this implication).$\blacksquare $
\end{Beweis}

\begin{theorem}
\label{cosep-counit}Let $(C,\Delta _{C},\varepsilon _{C})$ be a \emph{%
counital} $R$-coalgebra. If $_{R}C$ is locally projective, then the
following are equivalent:

\begin{enumerate}
\item $C$ is a coseparable $R$-coalgebra;

\item The coring $(\mathbf{C}:C^{\bullet })$ is left counital;

\item $C^{cop}$ is a coseparable $R$-coalgebra;

\item The coring $(\mathbf{C:C}^{\bullet })$ is right counital.
\end{enumerate}
\end{theorem}

\begin{Beweis}
Assume $_{R}C$ to be $R$-cogenerated.

$(1)\Rightarrow (2).$ Assume the $R$-coalgebra $C$ to be coseparable, so
that there exists a $(C,C)$-bicolinear map $\pi :C\otimes _{R}C\rightarrow C$
with $\pi \circ \Delta _{C}=\mathrm{id}_{C}.$ Consider the $R$-bilinear form%
\begin{equation*}
\gamma :C\times C\rightarrow R,\text{ }(c,d)\mapsto \varepsilon _{C}(\pi
(d\otimes _{R}c)),
\end{equation*}%
and the $R$-linear map%
\begin{equation}
\widetilde{\gamma }=\varepsilon _{C}\circ \pi \circ \tau :C\otimes
_{R}C\rightarrow R,\text{ }c\otimes _{R}d\mapsto <c,d>.
\end{equation}%
Then%
\begin{equation*}
\widetilde{\gamma }\circ \Delta _{C}^{\mathrm{tw}}=\varepsilon _{C}\circ \pi
\circ \tau \circ (\tau \circ \Delta _{C})=\varepsilon _{C}\circ \pi \circ
\Delta _{C}=\varepsilon _{C}\circ \mathrm{id}_{C}=\varepsilon _{C}.
\end{equation*}%
Moreover, we have%
\begin{equation*}
\begin{tabular}{lllll}
$(\mathrm{id}_{C}\otimes _{R}\widetilde{\gamma })\circ (\Delta _{C}^{\mathrm{%
tw}}\otimes _{R}\mathrm{id}_{C})$ & $=$ & $(\mathrm{id}_{C}\otimes
_{R}\varepsilon _{C})\circ (\mathrm{id}_{C}\otimes _{R}\pi \circ \tau )\circ
(\tau \circ \Delta _{C}\otimes _{R}\mathrm{id}_{C})$ &  &  \\
& $=$ & $(\mathrm{id}_{C}\otimes _{R}\varepsilon _{C})\circ \tau \circ (\pi
\otimes _{R}\mathrm{id}_{C})\circ (\mathrm{id}_{C}\otimes _{R}\Delta
_{C})\circ \tau $ &  &  \\
& $=$ & $(\mathrm{id}_{C}\otimes _{R}\varepsilon _{C})\circ \tau \circ
\Delta _{C}\circ \pi \circ \tau $ &  &  \\
& $=$ & $\pi \circ \tau $ &  &  \\
& $=$ & $(\mathrm{id}_{C}\otimes _{R}\varepsilon _{C})\circ \tau \circ
\Delta _{C}\circ \pi \circ \tau $ &  &  \\
& $=$ & $(\mathrm{id}_{C}\otimes _{R}\varepsilon _{C})\circ \tau \circ (%
\mathrm{id}_{C}\otimes _{R}\pi )\circ (\Delta _{C}\otimes _{R}\mathrm{id}%
_{C})\circ \tau $ &  &  \\
& $=$ & $(\varepsilon _{C}\otimes _{R}\mathrm{id}_{C})\circ (\pi \circ \tau
\otimes _{R}\mathrm{id}_{C})\circ (\mathrm{id}_{C}\otimes _{R}\tau \circ
\Delta _{C})$ &  &  \\
& $=$ & $(\widetilde{\gamma }\otimes _{R}\mathrm{id}_{C})\circ (\mathrm{id}%
_{C}\otimes _{R}\Delta _{C}^{\mathrm{tw}}).$ &  &
\end{tabular}%
\end{equation*}%
Hence%
\begin{equation*}
\widetilde{\gamma }\circ \Delta _{C}^{\mathrm{tw}}=\varepsilon _{C}\text{
and }(\mathrm{id}_{C}\otimes _{R}\widetilde{\gamma })\circ (\Delta _{C}^{%
\mathrm{tw}}\otimes _{R}\mathrm{id}_{C})=(\widetilde{\gamma }\otimes _{R}%
\mathrm{id}_{C})\circ (\mathrm{id}_{C}\otimes _{R}\Delta _{C}^{\mathrm{tw}}).
\end{equation*}%
By Lemma \ref{balanced}, the maps%
\begin{equation*}
\gamma ^{l}:\mathbf{C}\rightarrow C^{\bullet },\text{ }d\mapsto <-,d>\text{ (%
}\gamma ^{r}:\mathbf{C}\rightarrow C^{\bullet },\text{ }c\mapsto <c,->\text{)%
}
\end{equation*}%
are left (right) $C^{\bullet }$-linear. Moreover, we have for each $c\in C:$%
\begin{equation*}
\dsum \gamma ^{l}(c_{1})\rightharpoondown c_{2}=\dsum [\gamma
^{l}(c_{1})(c_{2})]c_{3}=\dsum <c_{2},c_{1}>c_{3}=\dsum \varepsilon
(c_{1})c_{2}=c.
\end{equation*}%
Consequently, $\gamma ^{l}:\mathbf{C}\rightarrow C^{\bullet }$ is a left
counit for the $C^{\bullet }$-coring $\mathbf{C}.$

$(2)\Rightarrow (3).$ Assume $\mathbf{C}$ is a $C^{\bullet }$-coring with
left counit $\varepsilon ^{l}$ and define%
\begin{equation*}
\widetilde{\pi }:C\otimes _{R}C\rightarrow C,\text{ }c\otimes _{R}d\mapsto
\varepsilon ^{l}(d)\rightharpoondown c.
\end{equation*}%
Then for all $f,g\in C^{\bullet },$ we have%
\begin{equation*}
\begin{tabular}{lllll}
$\widetilde{\pi }(f\rightharpoondown (c\otimes _{R}d)\leftharpoondown g)$ & $%
=$ & $\widetilde{\pi }((c\leftharpoondown g\otimes _{R}f\rightharpoondown
d)) $ & $=$ & $[\varepsilon ^{l}(f\rightharpoondown d)]\rightharpoondown
(c\leftharpoondown g)$ \\
& $=$ & $[f\bullet \varepsilon ^{l}(d)]\rightharpoondown (c\leftharpoondown
g)$ & $=$ & $([f\bullet \varepsilon ^{l}(d)]\rightharpoondown
c)\leftharpoondown g$ \\
& $=$ & $(f\rightharpoondown (\varepsilon ^{l}(d)\rightharpoondown
c))\leftharpoondown g$ & $=$ & $f\rightharpoondown \widetilde{\pi }(c\otimes
_{R}d)\leftharpoondown g.$%
\end{tabular}%
\end{equation*}%
So $\widetilde{\pi }:C\otimes _{R}C\rightarrow C$ is $(C^{\bullet
},C^{\bullet })$-bilinear, whence $(C,C)$-bicolinear by Proposition \ref{sg}%
. Moreover, we have for all $c\in C:$%
\begin{equation*}
(\widetilde{\pi }\circ \Delta _{C}^{\mathrm{tw}})(c)=\dsum \widetilde{\pi }%
(c_{2}\otimes _{R}c_{1})=\dsum \varepsilon ^{l}(c_{1})c_{2}=c.
\end{equation*}%
So, $C^{cop}$ is coseparable.

$(3)\Rightarrow (4).$ Analogous to the proof of $(1)\Rightarrow (2),$ we
conclude that the $C^{\ast }$-coring $\mathbf{C}^{cop}:=(C^{cop}:C^{\ast })$
is left counital, whence the coring $\mathbf{C}$ is right counital.

Analogous to the proof of $(2)\Rightarrow (3).\blacksquare $
\end{Beweis}

\textbf{Acknowledgement: }The author thanks King Fahd University of
Petroleum $\&\ $Minerals (KFUPM) for the financial support and the excellent
research facilities.

\end{document}